\documentclass[amsfonts,12pt]{article}

\usepackage{amsfonts,amsmath}
\usepackage{hyperref}

\newtheorem{theorem}{Theorem}
\newtheorem{corollary}[theorem]{Corollary}

\newtheorem{definition}[theorem]{Definition}

\DeclareMathOperator{\Div}{Div}
\DeclareMathOperator{\Ind}{Ind}

\def\ppp{{\mathbb{P}}}

\def\fff{{\mathbb{F}}}

\def\qqq{\mathbb{Q}}

\def\ccc{\mathbb{C}}

\def\pf{{\bf Proof}:\ }
\def\qed{$\Box$}
\def\div{{\rm div}}

\begin{document}

\author{David Joyner,
Amy Ksir\thanks{Mathematics Dept, USNA, Annapolis, MD 21402,
wdj@usna.edu and ksir@usna.edu},
Roger Vogeler\thanks{Mathematics Dept.,
Ohio State Univ., vogeler@math.ohio-state.edu}
}
\title{Group representations on \\ Riemann-Roch spaces \\
of some Hurwitz curves }
\date{11-14-2006}

\maketitle

\begin{abstract}
Let $q>1$ denote an integer relatively prime to $2,3,7$ and for
which $G=PSL(2,q)$ is a Hurwitz group for a smooth projective curve
$X$ defined over $\ccc$. We compute the $G$-module structure of the
Riemann-Roch space $L(D)$, where $D$ is an invariant divisor on $X$
of positive degree. This depends on a computation of the
ramification module, which we give explicitly. In particular, we
obtain the decomposition of $H^1(X,\ccc)$ as a $G$-module.

\end{abstract}

\vskip .5in

\tableofcontents

\vskip .3in

\section{Introduction}

Let $X$ be a smooth projective curve over an algebraically closed
field $k$, and let $k(X)$ denote the function field of $X$ (the
field of rational functions on $X$).  If $D$ is any divisor on $X$
then the Riemann-Roch space $L(D)$ is a finite dimensional
$k$-vector space given by

\[
L(D)=L_X(D)= \{f\in k(X)^\times \ |\ \div(f)+D\geq 0\}\cup \{0\},
\]
where $\div(f)$ denotes the (principal) divisor of the function
$f\in k(X)$.   If $G$ is a finite group of automorphisms of $X$,
then $G$ has a natural action on $k(X)$, and on the group $\Div(X)$
of divisors on $X$.  If $D$ is a $G$-invariant divisor, then $G$
also acts on the vector space $L(D)$, making it into a
$k[G]$-module.

The problem of finding the $k[G]$-module structure of $L(D)$ was
first considered in the case where $k=\ccc$ and $D$ is canonical,
i.e. $L(D)$ is the space of holomorphic differentials on $X$.  This
problem was solved by Hurwitz for $G$ cyclic, and then by Chevalley
and Weil for general $G$.  More generally, the problem has been
solved by work of Ellingsrud and L{\o}nsted \cite{EL}, Kani
\cite{K}, Nakajima \cite{N}, and Borne \cite{B}. This has resulted
in the following equivariant Riemann-Roch formula for the class of
$L(D)$ (denoted by square brackets) in the Grothendieck group
$R_k(G)$, in the case where $D$ is non-special:
\begin{equation}
\label{eqn:Borne}
[L(D)]=(1-g_{X/G})[k[G]]+[\deg_{eq}(D)]-[\tilde{\Gamma}_G].
\end{equation}
Here $g_{X/G}$ is the genus of $X/G$, $\deg_{eq}(D)$ is the
equivariant degree of $D$, and $\tilde{\Gamma}_G$ is the (reduced)
ramification module  (this notation will be defined in sections
\ref{sec:rammod} and \ref{sec:equivdeg}).

Explicitly computing the $k[G]$-module structure of $L(D)$ in
specific cases is of interest currently due to advances in the
theory of algebraic-geometric codes.  Permutation decoding
algorithms use this information to increase their efficiency.

In this paper, we consider the case where $X$ is a Hurwitz curve
with automorphism group $G=PSL(2,q)$ for some prime power $q$, over
$k=\ccc$.  Using the equivariant Riemann-Roch formula above
(\ref{eqn:Borne}) and the representation theory of $PSL(2,q)$, we
compute explicitly the $\ccc[G]$-module structure of $L(D)$ for a
general invariant effective divisor $D$.  In the case where $D$ is a
canonical divisor, this yields an explicit computation for the
$\ccc[G]$-module structure of $H^{1}(X,\ccc)$.

We are also interested in rationality questions.  We find that
$\tilde{\Gamma}_G$ has a $\qqq[G]$-module structure, and therefore
may be computed more simply (see Joyner and Ksir \cite{JK1}), as
follows:
\begin{equation}
\label{eqn:JKrammod}
\tilde{\Gamma}_G
=\bigoplus_{\pi\in G^*} \left[ \sum_{\ell=1}^L
({\rm dim}\, \pi -{\rm dim}\,
(\pi^{H_\ell})) \frac{R_\ell}{2}\right]\pi .
\end{equation}
The sum is over all conjugacy classes of cyclic subgroups of $G$,
$H_\ell$ is a representative cyclic subgroup, $\pi^{H_\ell}$
indicates the fixed part of $\pi$ under the action of $H_\ell$, and
$R_\ell$ denotes the number of branch points in $Y$ over which the
decomposition group is conjugate to $H_\ell$.  For some but not all
divisors $D$, $L(D)$ has a $\qqq[G]$-module structure, and may also
be computed more simply.

The organization of this paper is as follows.  In section 2, we
recall some facts about Hurwitz curves and Hurwitz groups.  In
section 3, we review the representation theory of $PSL(2,q)$, and
compute the induced characters necessary for the following section.
Our main results are in section $4$, where we compute the
ramification module, the equivariant degree for any invariant
divisor $D$, and thus the structure of $L(D)$.  At the end of
section 4 we compute the $\ccc[G]$-module structure of
$H^1(X,\ccc)$.  In section 5, we discuss rationality questions,
using the results of \cite{JK1} to give more streamlined formulas
for the ramification module, and in some cases for $L(D)$.

\section{Hurwitz curves}
\label{sec:hurwitz}

The automorphism group $G$ of a smooth projective
curve of genus $g>1$ over an
algebraically closed field $k$ of characteristic zero
satisfies the {\it Hurwitz bound}

\[
|G|\leq 84\cdot (g-1).
\]
A curve which attains this bound is called a
{\it Hurwitz curve} and its automorphism group
is called a {\it Hurwitz group}.

\subsection{Classification}

The number of distinct Hurwitz groups is infinite, and
to each one corresponds a finite number of Hurwitz curves.
Nevertheless, these curves are quite rare;
in particular, the Hurwitz genus values are known to form
a rather sparse set of positive integers (see Larsen \cite{L}).

Hurwitz groups are precisely those groups which occur
as non-trivial finite homomorphic images of the 2,3,7-triangle group
\[
\Delta =\langle a,b : a^2=b^3=(ab)^7=1\rangle.
\]
This is most naturally viewed as the group of orientation-preserving
symmetries of the tiling of the hyperbolic plane
$\bold{H}$ generated by reflections in the sides of a fundamental
triangle having angles $\pi/2$, $\pi/3$, and $\pi/7$.
Each proper normal finite-index subgroup $K\triangleleft\Delta$
corresponds to a Hurwitz group $G=\Delta/K$. The associated
Hurwitz curve now appears (with $k=\ccc$) as a compact hyperbolic
surface $\bold{H}/K$ regularly tiled by a finite number
of copies of the fundamental triangle. $G$ is the group of
orientation-preserving symmetries of this tiling,
with fundamental domain consisting of one fundamental
triangle plus one reflected triangle. (From this perspective,
the Hurwitz bound simply says that there is no smaller
polygon which gives a regular tiling of $\bold{H}$.)

We note that $\Delta$ has only a small number of torsion
elements (up to conjugacy). These are the non-trivial
powers of $a$, $b$, and $ab$. Each acts as a rotation of
order 2, 3, or 7, and has as its fixed point one vertex of
(some copy of) the fundamental triangle. Clearly no other
point of the tiling can occur as a fixed point; this is
true both for the tiling of $\bold{H}$ and the induced
tilings on the quotient surfaces. In other words, all
points {\it other} than the tiling vertices have trivial stabilizer.

It follows easily from the above presentation for $\Delta$
that a group is Hurwitz if and only if it is generated
by two elements having orders 2 and 3, and whose
product has order 7. This characterization has made
possible much of the work in classifying Hurwitz
groups. The most relevant for our investigation is the
following result of Macbeath (see \cite{M}):

\begin{verse}
    The simple group $PSL(2,q)$ is Hurwitz in exactly three cases:\\
    \quad i) $q=7$;\\
    \quad ii) $q$ is prime, with $q\equiv \pm 1 \pmod{7}$;\\
    \quad iii) $q=p^3$, with $p$ prime and $p\equiv \pm 2,\pm 3$ (mod 7).
\end{verse}
In particular, $PSL(2,8)$ and $PSL(2,27)$ are Hurwitz groups. We
shall require that $q$ be relatively prime to $2\cdot 3\cdot 7$, but
this excludes just three possibilities, namely $q\in \{7,8,27\}$.
Note that in all of the cases we consider, $q \equiv \pm 1
\pmod{7}$.

The order of $PSL(2,q)$ (for odd $q$) is $q(q^2-1)/2$. Hence we obtain
\[
g=1+\frac{q(q^2-1)}{168}
\]
as the genus of the corresponding curve(s).

For completeness, we remark that there are three distinct Hurwitz curves when $q$ is prime (apart from $q=7$), and just one when $q=p^3$. However, this has no bearing on the representations that we study.

In addition, there are other known families of Hurwitz groups. For example, all Ree groups are Hurwitz, as are all but finitely many of the alternating groups. See Conder \cite{C} for a summary of such results.

\subsection{Ramification data}
\label{sec:ram}

Let $X$ be a Hurwitz curve with automorphism group $G$ and
let

\begin{equation}
\label{eqn:psi}
\psi:X\rightarrow Y=X/G
\end{equation}
denote the quotient map. By again viewing $X$ as a hyperbolic
surface, the ramification data are easily deduced. The quotient $Y$
is formed by one fundamental triangle and its mirror image, with the
natural identifications on their boundaries. Hence it is a surface
of genus 0 with 3 metric singularities. Thus $\psi$ has exactly
three branch points.  The stabilizer subgroups of the corresponding
ramification points in $X$ are cyclic, of orders $2$, $3$, and $7$.
We label the three branch points $P_1$, $P_2$, and $P_3$, so that if
$P \in \psi^{-1}(P_1)$, then $P$ has stabilizer subgroup of order
$2$, if $P \in \psi^{-1}(P_2)$, $P$ has stabilizer subgroup of order
$3$, and if $P \in \psi^{-1}(P_3)$, $P$ has stabilizer subgroup of
order $7$.

\section{Representation theory of $PSL(2,q)$}
\label{sec:rep thy}

\subsection{General theory on representations of PSL(2,q)}

We first review the representation theory of $G=PSL(2,q)$ over
$\ccc$, following the treatment in \cite{FH}, to fix notation.

Let $\fff=GF(q)$ be the field with $q$ elements. The group
$PSL(2,q)$ has $3+(q-1)/2$ conjugacy classes of elements. Let
$\varepsilon \in {\fff}$ be a generator for the cyclic group
${\fff}^{\times}$.  Then each conjugacy class will have a
representative of exactly one of the following forms:

\begin{equation}
\label{eqn:conjclasses}
 \left(
\begin{array}{cc}
1 & 0\\
0 &1
\end{array}
\right),\
\left(
\begin{array}{cc}
x & 0\\
0 & x^{-1}
\end{array}
\right),\
\left(
\begin{array}{cc}
1 & 1\\
0 & 1
\end{array}
\right),\
\left(
\begin{array}{cc}
1 & \varepsilon\\
0 & 1
\end{array}
\right),\
\left(
\begin{array}{cc}
x & \varepsilon y\\
y & x
\end{array}
\right).
\end{equation}

The irreducible representations of $PSL(2,q)$ include the trivial
representation $\mathbf{1}$ and one irreducible $V$ of dimension
$q$. All but two of the others fall into two types:  representations
$W_{\alpha}$ of dimension $q+1$ (``principal series''), and
$X_{\beta}$ of dimension $q-1$ (``discrete series''). The principal
series representations $W_{\alpha}$ are indexed by homomorphisms
$\alpha: {\fff}^{\times} \to \ccc^{\times}$ with $\alpha(-1)=1$. The
discrete series representations $X_{\beta}$ are indexed by
homomorphisms $\beta: T \to \ccc^{\times}$ with $\beta(-1)=1$, where
$T$ is a cyclic subgroup of order $q+1$ of
${\fff}(\sqrt{\varepsilon})^{\times}$. The characters of these are
as follows:

{\footnotesize{
\[
\begin{array}[ht]{r||c|c|cc|c}
& \left( \begin{array}{cc} 1& 0 \\ 0 & 1 \end{array} \right) &
\left( \begin{array}{cc} x& 0 \\ 0 & x^{-1} \end{array} \right)&
\left( \begin{array}{cc} 1& 1 \\ 0 & 1 \end{array} \right)&
\left( \begin{array}{cc} 1& \varepsilon \\ 0 & 1 \end{array} \right)&
\left( \begin{array}{cc} x& \varepsilon y \\ y & x \end{array} \right)\\
\hline \hline
\mathbf{1} & 1 & 1 & 1 & 1 &1 \\
\hline X_{\beta} & q-1 & 0 & -1 & -1 &
-\beta(x+\sqrt{\varepsilon}y)-\beta(x-\sqrt{\varepsilon}y) \\
\hline V & q & 1 & 0 & 0 & -1 \\
\hline W_{\alpha} & q+1 & \alpha(x) + \alpha(x^{-1}) & 1 & 1 & 0 \\
\end{array}
\]
}}

Let $\zeta$ be a primitive
$q$th root of unity in $\ccc$.  Let $\xi$ and $\xi'$ be defined by

\begin{equation}
\label{eqn:qq'}
\xi = \sum_{\left(\frac{a}{q}\right)=1} \zeta^a \mbox{
and } \xi' = \sum_{\left(\frac{a}{q}\right)=-1} \zeta^a,
\end{equation}
where the sums are over the quadratic residues and nonresidues
$\pmod q$, respectively. If $q \equiv 1$ mod 4, then the principal
series representation $W_{\alpha_0}$ corresponding to

\[
\begin{array}{ccc}
\alpha_0:{\fff}^{\times} & \to & \ccc^{\times} \\
\varepsilon & \mapsto & -1
\end{array}
\]
is not irreducible, but splits into two irreducibles $W'$ and $W''$,
each of dimension $(q+1)/2$.
Their characters satisfy:

\[
\begin{array}[ht]{r||c|c|cc|c}
& \left( \begin{array}{cc} 1& 0 \\ 0 & 1 \end{array} \right) &
\left( \begin{array}{cc} x& 0 \\ 0 & x^{-1} \end{array} \right)&
\left( \begin{array}{cc} 1& 1 \\ 0 & 1 \end{array} \right)& \left(
\begin{array}{cc} 1& \varepsilon \\ 0 & 1 \end{array} \right)&
\left( \begin{array}{cc} x& \varepsilon y \\ y & x \end{array} \right)\\
\hline \hline W' & \frac{q+1}{2} & \alpha_0(x) & 1+\xi
& 1+\xi' & 0 \\
\hline W'' & \frac{q+1}{2} & \alpha_0(x) & 1+\xi'
& 1+\xi & 0 \\
\end{array}
\]

Let $\tau$ denote a generator of $T$. Similarly, if $q \equiv 3$ mod
4, then the discrete series representation $X_{\beta_0}$
corresponding to

\[
\begin{array}{ccc}
\beta_0:T & \to & \ccc^{\times} \\
\tau & \mapsto & -1
\end{array}
\]
splits into two irreducibles $X'$ and $X''$, each of dimension
$(q-1)/2$. Their characters satisfy:

\[
\begin{array}[ht]{r||c|c|cc|c}
& \left( \begin{array}{cc} 1& 0 \\ 0 & 1 \end{array} \right) &
\left( \begin{array}{cc} x& 0 \\ 0 & x^{-1} \end{array} \right)&
\left( \begin{array}{cc} 1& 1 \\ 0 & 1 \end{array} \right)& \left(
\begin{array}{cc} 1& \varepsilon \\ 0 & 1 \end{array} \right)&
\left( \begin{array}{cc} x& \varepsilon y \\ y & x \end{array} \right)\\
\hline \hline X' & \frac{q-1}{2} & 0 & \xi
& \xi' & -\beta_0(x+y\sqrt{\varepsilon}) \\
\hline X'' & \frac{q-1}{2} & 0 & \xi'
& \xi & -\beta_0(x+y\sqrt{\varepsilon}) \\
\end{array}
\]
According to Janusz \cite{Ja}, the Schur index of each irreducible
representation of $G$ is $1$.

There is a ``Galois action'' on the set of equivalence classes of
irreducible representations of $G$ as follows. Let $\chi$ denote an
irreducible character. The character values $\chi(g)$ lie in
$\qqq(\mu)$, where $\mu$ is a primitive $m^{th}$ root of unity and
$m=q(q^2-1)/4$.  Let ${\cal{G}}=Gal(\qqq(\mu)/\qqq)$ denote the
Galois group.  For each integer $j$ relatively prime to $m$, there
is an element $\sigma_j$ of ${\cal{G}}$ taking $\mu$ to $\mu^j$.
This Galois group element will act on representations by taking a
representation with character values $(a_1, \ldots, a_n)$ to a
representation with character values $(\sigma_j(a_1), \ldots,
\sigma_j(a_n))$. Representations with rational character values will
be fixed under this action.  Because the Schur index of each
representation is $1$, representations with rational character
values will be defined over $\qqq$.

The action of the Galois group $\mathcal{G}$ can easily be seen from
the character table. It will fix the trivial representation and the
$q$-dimensional representation $V$. Its action permutes the set of
$q-1$-dimensional ``principal series'' representations $X_{\beta}$,
and the set of $q+1$-dimensional ``discrete series'' representations
$W_{\alpha}$. In the case $q \equiv 1\pmod 4$, the Galois group will
exchange the two $(q+1)/2$-dimensional representations $W'$ and
$W''$; if $q \equiv 3 \pmod 4$, the Galois group will exchange the
two $(q-1)/2$-dimensional representations $X'$ and $X''$.

\subsection{Induced characters}
\label{sec:indchars}

We will be interested in the induced characters from subgroups of
orders $2$, $3$, and $7$.  For each value of $q$, each of these
subgroups is unique up to conjugacy; we can choose subgroups $H_2$
of order $2$, $H_3$ of order $3$, and $H_7$ of order $7$ that are
generated by elements of the form
\begin{equation*}
\left(
\begin{array}{cc}
x & 0\\
0 & x^{-1}
\end{array}
\right)
\mbox{ or }
\left(
\begin{array}{cc}
x & \varepsilon y\\
y & x
\end{array}
\right).
\end{equation*}
Which of these two forms each generator will take depends on $q$ mod
4, mod 3, and mod 7, respectively. Recall that we defined generators
$\varepsilon$ of the cyclic group $\fff^{\times}$, of order $q-1$,
and $\tau$ of the cyclic group $T \subseteq
\fff(\sqrt{\varepsilon})^{\times}$ of order $q+1$, respectively.  We
define numbers $i$, $\omega$, and $\phi$ to be primitive roots of
unity as follows.

When $q \equiv 1 \pmod 4$, let $i$ denote an element in
${\fff}^{\times}$ whose square is $-1$ (one can take
$i=\varepsilon^{(q-1)/4}$). Then the subgroup $H_2$ of order $2$ in
$PSL(2,q)$ is generated by
\[
\left(
\begin{array}{cc}
i & 0\\
0 & i^{-1}
\end{array}
\right ).
\]
If $q \equiv 3 \pmod 4$, then we take $i=x_i+\sqrt{\varepsilon}y_i$
to be an element of $T$ whose square is $-1$ (one can take
$i=\tau^{(q+1)/4}$).  Then the subgroup $H_2$ of order $2$ in
$PSL(2,q)$ is generated by
\[
\left(
\begin{array}{cc}
x_i & \varepsilon y_i\\
y_i & x_i
\end{array}
\right ).
\]

Similarly, we define $\omega$ to be a primitive $6$th root of unity.
In the case where $q \equiv 1 \pmod 3$, we can take
$\omega=\varepsilon^{(q-1)/6} \in \fff^{\times}$.  When $q \equiv -1
\pmod 3$, we take $\omega = x_{\omega} + \sqrt{\varepsilon}
y_{\omega} = \tau^{(q+1)/6} \in T$. The subgroup $H_3$ of order $3$
in $PSL(2,q)$ will then be generated by
\[
\left(
\begin{array}{cc}
\omega & 0\\
0 & \omega^{-1}
\end{array}
\right ),\mbox{ if } \  q \equiv 1 \pmod 3, \mbox{ or } \left(
\begin{array}{cc}
x_{\omega} & \varepsilon y_{\omega}\\
y_{\omega} & x_{\omega}
\end{array}
\right ),\mbox{ if } \  q \equiv -1 \pmod 3.
\]

Lastly, we want to define $\phi$ to be a primitive $14$th root of
unity. Recall that $q \equiv \pm 1 \pmod 7$.  If $q \equiv 1 \pmod
7$, then we can take $\phi = \varepsilon^{(q-1)/14} \in
\fff^{\times}$, and if $q \equiv -1 \pmod 7$, then we can take $\phi
= x_{\phi} + \sqrt{\varepsilon} y_{\phi} = \tau^{(q+1)/14} \in T$.
The subgroup $H_7$ of order $7$ in $PSL(2,q)$ will then be generated
by
\[
\left(
\begin{array}{cc}
\phi & 0\\
0 & \phi^{-1}
\end{array}
\right ), \ q \equiv 1 \pmod 3, \mbox{ or } \left(
\begin{array}{cc}
x_{\phi} & \varepsilon y_{\phi}\\
y_{\phi} & x_{\phi}
\end{array}
\right ), \ q \equiv -1 \pmod 3.
\]

With these definitions, it is easy to compute the restrictions of
the irreducible representations of $PSL(2,q)$ to the subgroups
above. We omit the details, but the computations for the groups of
order $2$ and $3$ are given in \cite{JK2}, and the computation for
the group of order $7$ is very similar.  Using Frobenius
reciprocity, we then obtain the corresponding induced
representations. In each case, we denote a primitive character of
the cyclic group $H_k$ by $\theta_k$.

\subsubsection{Induced characters from $H_2$}
The induced representations from the nontrivial character of $H_2$
are given below. The multiplicities depend on $q \pmod 8$.  Note
that most representation have the same multiplicity as $V$. When $i
\in \fff^{\times}$, i.e. when $q \equiv 1 \pmod 4$, the multiplicity
of a discrete series representation $W_{\alpha}$ depends on the sign
of $\alpha(i)$.  Recall that $\alpha(-1)=1$, so $\alpha(i) = \pm 1$.
the multiplicity of $W_{\alpha}$ will be the same as the
multiplicity of $V$ if $\alpha(i)=1$ and one larger if $\alpha(i) =
-1$. Similarly, when $q \equiv 3 \pmod 4$ and $i \in T$, the
multiplicity of a principal series representation $X_{\beta}$
depends on the sign of $\beta(i)$.  In this case the multiplicity of
$X_{\beta}$ will be the same as the multiplicity of $V$ when
$\beta(i)=1$, and one less if $\beta(i) = -1$.  Lastly, the signs of
$\alpha_0(i)$ or $\beta_0(i)$ depend on $q \pmod 8$ and determine
the multiplicities of $W'$ and $W''$ or $X'$ and $X''$,
respectively.  A similar pattern will hold for the induced
representations from $H_3$ and $H_7$.

For $q \equiv 1 \pmod 8$, {\footnotesize{
\begin{eqnarray*}
Ind_{H_2}^{G} \theta_2 & = & \frac{q-1}{2} \left [ \frac{1}{2}(W' +
W'') + \sum_{\beta} X_{\beta} +  V + \sum_{\alpha(i)=1} W_{\alpha}
\right ] + \frac{q+3}{2} \sum_{\alpha(i)=-1} W_{\alpha}.
\end{eqnarray*}
}} For $q \equiv 3 \pmod 8$, {\footnotesize{
\begin{eqnarray*}
Ind_{H_2}^{G} \theta_2 & = & \frac{q+1}{2} \left [ \sum_{\beta(i)=1}
X_{\beta} + V + \sum_{\alpha} W_{\alpha} \right ] + \frac{q-3}{2}
\left [ \frac{1}{2}(X' + X'') + \sum_{\beta(i)=-1} X_{\beta} \right
].
\end{eqnarray*}
}} For $q \equiv 5 \pmod 8$, {\footnotesize{
\begin{eqnarray*}
Ind_{H_2}^{G} \theta_2 & = & \frac{q-1}{2} \left [ \sum_{\beta}
X_{\beta} + V + \sum_{\alpha(i)=1} W_{\alpha} \right ] +
\frac{q+3}{2} \left [ \frac{1}{2}(W' + W'') + \sum_{\alpha(i)=-1}
W_{\alpha} \right ].
\end{eqnarray*}
}} And for $q \equiv 7 \pmod 8$, {\footnotesize{
\begin{eqnarray*}
Ind_{H_2}^{G} \theta_2 & = & \frac{q+1}{2} \left [ \frac{1}{2} (X' +
X'') + \sum_{\beta(i)=1} X_{\beta} + V + \sum_{\alpha} W_{\alpha}
\right ] + \frac{q-3}{2} \sum_{\beta(i)=-1} X_{\beta}.
\end{eqnarray*} }}

\subsubsection{Induced characters from $H_3$}
The induced representations from the two nontrivial characters
$\theta_3$ and $\theta_3^2$ of $H_3$ are the same.  In this case the
multiplicities depend on $q \pmod {12}$, which determines whether
the $6$th root of unity $\omega$ is in $\fff^{\times}$, or in $T
\subset \fff(\sqrt{\varepsilon})^{\times}$.  Now the multiplicity of
a discrete (resp. principal) series representation $W_{\alpha}$
(resp. $X_{\beta}$) will be the same as the multiplicity of $V$ if
$\alpha(\phi)=1$ (resp. $\beta(\phi)=1$) and one larger (resp.
smaller) if $\alpha(\phi) = e^{\frac{\pm 2 \pi i}{3}}$ (resp.
$\beta(\phi) = e^{\frac{\pm 2 \pi i}{3}}$).  The signs of
$\alpha_0(\omega)$ or $\beta_0(\omega)$ depend on $q \pmod {12}$ and
determine the multiplicities of $W'$ and $W''$ or $X'$ and $X''$,
respectively.

If $q \equiv 1 \pmod {12}$, we have

{\footnotesize{
\begin{equation*}
Ind_{H_3}^{G} \theta_3 = \frac{q-1}{3} \left [ \frac{1}{2}(W' + W'')
+ \sum_{\beta} X_{\beta} + V + \sum_{\alpha(\omega)=1} W_{\alpha}
 \right ] + \frac{q+2}{3} \sum_{\alpha(\omega)
= e^{\frac{\pm 2 \pi i}{3}}} W_{\alpha}.
\end{equation*} }}

If $q \equiv 5 \pmod {12}$, we have
{\footnotesize{
\begin{equation*}
Ind_{H_3}^{G} \theta_3  =  \frac{q+1}{3} \left [ \frac{1}{2} (W' +
W'') + \sum_{\beta(\omega)=1} X_{\beta} + V + \sum_{\alpha}
W_{\alpha} \right ] + \frac{q-2}{3} \sum_{\beta(\omega) = 1}
X_{\beta}.
\end{equation*} }}

If $q \equiv 7 \pmod {12}$, we have
{\footnotesize{
\begin{equation*}
Ind_{H_3}^{G} \theta_3 = \frac{q-1}{3} \left [ \frac{1}{2}(X' + X'')
+ \sum_{\beta} X_{\beta} + V + \sum_{\alpha(\omega)=1} W_{\alpha}
\right ] + \frac{q+2}{3} \sum_{\alpha(\omega) = e^{\frac{\pm 2 \pi
i}{3}}} W_{\alpha}.
\end{equation*} }}

And if $q \equiv 11 \pmod {12}$, we have {\footnotesize{
\begin{equation*}
Ind_{H_3}^{G} \theta_3  =  \frac{q+1}{3} \left [\frac{1}{2} (X' +
X'') + \sum_{\beta(\omega)=1} X_{\beta} + V + \sum_{\alpha}
W_{\alpha}
 \right ] + \frac{q-2}{3} \sum_{\beta(\omega)=
e^{\frac{\pm 2 \pi i}{3}}} X_{\beta}.
\end{equation*} }}

\subsubsection{Induced characters from $H_7$}
For $H_7$, the induced representations from the six nontrivial
characters $\theta_7^k$ are not all the same, but depend on $k$.
These representations also depend on $q \pmod {28}$, which
determines whether the $14$th root of unity $\phi$ is in
$\fff^{\times}$ or $\fff(\sqrt{\varepsilon})^{\times}$. For an
induced nontrivial character $\Ind_{H_7}^{G} \theta_7^k$, the
multiplicity of a discrete (resp. principal) series representation
$W_{\alpha}$ (resp. $X_{\beta}$) will be the same as the
multiplicity of $V$ if $\alpha(\phi) \neq e^{\pm \frac{2 \pi i
k}{7}}$ (resp. $\beta(\phi) \neq e^{\pm \frac{2 \pi i k}{7}}$) and
one larger (resp. smaller) if $\alpha(\phi) = e^{\pm \frac{2 \pi i
k}{7}}$ (resp. $\beta(\phi) = e^{\pm \frac{2 \pi i k}{7}}$). The
signs of $\alpha_0(\phi)$ or $\beta_0(\phi)$ depend on $q \pmod
{28}$ and determine the multiplicities of $W'$ and $W''$ or $X'$ and
$X''$, respectively.

If $q \equiv 1 \pmod {28}$, we have

{\footnotesize{
\begin{equation*}
Ind_{H_7}^{G} \theta_7^k = \frac{q-1}{7} \left [\frac{1}{2} (W' +
W'') + \sum_{\beta} X_{\beta} + V + \sum_{\alpha(\phi) \neq e^{\pm
\frac{2 \pi i k}{7}}} W_{\alpha} \right ] + \frac{q+6}{7}
\sum_{\alpha(\phi) = e^{\pm \frac{2 \pi i k }{7}}} W_{\alpha}.
\end{equation*} }}

If $q \equiv 13 \pmod {28}$, we have

{\footnotesize{
\begin{equation*}
Ind_{H_7}^{G} \theta_7^k = \frac{q+1}{7} \left [ \frac{1}{2} (W' +
W'') + \sum_{\beta(\phi) \neq e^{\pm \frac{2 \pi i k}{7}}} X_{\beta}
+ V + \sum_{\alpha} W_{\alpha} \right ] + \frac{q-6}{7}
\sum_{\beta(\phi) = e^{\pm \frac{2 \pi i k}{7}}} X_{\beta}.
\end{equation*} }}

If $q \equiv 15 \pmod {28}$, we have

{\footnotesize{
\begin{equation*}
Ind_{H_7}^{G} \theta_7^k = \frac{q-1}{7} \left [ \frac{1}{2} (X' +
X'') + \sum_{\beta} X_{\beta} + V + \sum_{\alpha(\phi) \neq e^{\pm
\frac{2 \pi i k}{7}}} W_{\alpha} \right ] + \frac{q+6}{7}
\sum_{\alpha(\phi) = e^{\pm \frac{2 \pi i k }{7}}} W_{\alpha}.
\end{equation*} }}

And if $q \equiv 27 \pmod {28}$, we have

{\footnotesize{
\begin{equation*}
Ind_{H_7}^{G} \theta_7^k = \frac{q+1}{7} \left [ \frac{1}{2} (X' +
X'') + \sum_{\beta(\phi) \neq e^{\pm \frac{2 \pi i k}{7}}} X_{\beta}
+ V + \sum_{\alpha} W_{\alpha} \right ] + \frac{q-6}{7}
\sum_{\beta(\phi) = e^{\pm \frac{2 \pi i k}{7}}} X_{\beta}.
\end{equation*} }}

\section{The Riemann-Roch space as a $G$-module}

Now we have all of the pieces we need to compute the $G$-module
structure of the Riemann-Roch space $L(D)$ of a general
$G$-invariant divisor $D$.  We will first compute the ramification
module, which does not depend on $D$.  We will then compute the
equivariant degree of $D$, and use the equivariant Riemann-Roch
formula (\ref{eqn:Borne}) to compute $L(D)$.

\subsection{Ramification module}
\label{sec:rammod}

The ramification module introduced by Kani \cite{K} and Nakajima
\cite{N} is defined by
\[
\Gamma_G=\sum_{P\in X_{\rm{ram}}} \Ind_{G_P}^G\left( \sum_{\ell
=1}^{e_P-1} \ell\theta_P^\ell \right),
\]
where the first sum is over the ramification points of
$\psi:X\rightarrow Y=X/G$, and $\theta_P$ is the ramification
character at a point $P$. Both Kani and Nakajima showed that there
is a $G$-module $\tilde{\Gamma}_G$ such that $\Gamma_G \simeq
\bigoplus_{|G|} \tilde{\Gamma}_G$.  Because $\Gamma_G$ does not
figure in our calculations, we abuse notation and refer to
$\tilde{\Gamma}_G$ as the {\it ramification module}.

Recall from section \ref{sec:ram} that $\psi:X\rightarrow Y=X/G$ has
three branch points, $P_1$, $P_2$, and $P_3$.  If $P \in
\psi^{-1}(P_1)$, $G_P$ has order $2$, so there are $\frac{|G|}{2}$
ramification points where $G_P$ is conjugate to $H_2$. If $P \in
\psi^{-1}(P_2)$, $G_P$ has order $3$, so there are $\frac{|G|}{3}$
ramification points where $G_P$ is conjugate to $H_3$, and if $P \in
\psi^{-1}(P_3)$, $G_P$ has order $7$, so there are $\frac{|G|}{7}$
ramification points where $G_P$ is conjugate to $H_7$.  Thus

\begin{equation}
\tilde{\Gamma}_G= \frac{1}{|G|} \left (\frac{|G|}{2}\Ind_{H_2}^{G}
\theta_2 +\frac{|G|}{3} \sum_{\ell =1}^{2} \ell
\Ind_{H_3}^{G}\theta_3^{\ell}  +\frac{|G|}{7} \sum_{\ell=1}^{6} \ell
\Ind_{H_7}^{G} \theta_7^{\ell} \right ).
\end{equation}

To compute this, we break it into three pieces:

\begin{eqnarray*}
\tilde{\Gamma}_G & = & \Gamma_{H_2} + \Gamma_{H_3} + \Gamma_{H_7}, \\
\Gamma_{H_2} & = & \frac{1}{2} \Ind_{H_2}^{G} \theta_2, \\
\Gamma_{H_3} & = & \frac{1}{3} ( \Ind_{H_3}^{G} \theta_3 + 2
\Ind_{H_3}^{G} \theta_3^2), \\
\Gamma_{H_7} & = & \frac{1}{7} ( \Ind_{H_7}^{G} \theta_7 + 2
\Ind_{H_7}^{G} \theta_7^2 + 3 \Ind_{H_7}^{G} \theta_7^3 \\
& & + 4 \Ind_{H_7}^{G} \theta_7^4 + 5 \Ind_{H_7}^{G} \theta_7^5 + 6
\Ind_{H_7}^{G} \theta_7^6 ).
\end{eqnarray*}

Each piece is then computed from the induced characters in section
\ref{sec:indchars}.  $\Gamma_{H_2}$ depends on $q \pmod 8$.

For $q \equiv 1 \pmod 8$, {\footnotesize{
\begin{eqnarray*}
\Gamma_{H_2} & = & \frac{q-1}{4} \left [ \frac{1}{2}(W' + W'') +
\sum_{\beta} X_{\beta} +  V + \sum_{\alpha(i)=1} W_{\alpha} \right ]
+ \frac{q+3}{4} \sum_{\alpha(i)=-1} W_{\alpha}.
\end{eqnarray*}
}} For $q \equiv 3 \pmod 8$, {\footnotesize{
\begin{eqnarray*}
\Gamma_{H_2} & = & \frac{q+1}{4} \left [ \sum_{\beta(i)=1} X_{\beta}
+ V + \sum_{\alpha} W_{\alpha} \right ] + \frac{q-3}{4} \left [
\frac{1}{2}(X' + X'') + \sum_{\beta(i)=-1} X_{\beta} \right ].
\end{eqnarray*}
}} For $q \equiv 5 \pmod 8$, {\footnotesize{
\begin{eqnarray*}
\Gamma_{H_2} & = & \frac{q-1}{4} \left [ \sum_{\beta} X_{\beta} + V
+ \sum_{\alpha(i)=1} W_{\alpha} \right ] + \frac{q+3}{4} \left [
\frac{1}{2}(W' + W'') + \sum_{\alpha(i)=-1} W_{\alpha} \right ].
\end{eqnarray*}
}} And for $q \equiv 7 \pmod 8$, {\footnotesize{
\begin{eqnarray*}
\Gamma_{H_2} & = & \frac{q+1}{4} \left [ \frac{1}{2} (X' + X'') +
\sum_{\beta(i)=1} X_{\beta} + V + \sum_{\alpha} W_{\alpha} \right ]
+ \frac{q-3}{4} \sum_{\beta(i)=-1} X_{\beta}.
\end{eqnarray*} }}

The contribution $\Gamma_{H_3}$ of $H_3$ to the ramification module
is
\begin{equation*}
\Gamma_{H_3} = \frac{1}{3} \left ( \Ind_{H_3}^{G} \theta_3 + 2
\Ind_{H_3}^{G} \theta_3^2 \right) = \Ind_{H_3}^{G} \theta_3,
\end{equation*}
since $\Ind_{H_3}^{G} \theta_3$ and $\Ind_{H_3}^{G} \theta_3^2$ are
the same.  This character was computed in section
\ref{sec:indchars}.

For $H_7$, the induced representations from the six nontrivial
characters $\theta_7^k$ are not all the same.  However, the
representations $\Ind_{H_7}^{G} \theta_7^k$ and $\Ind_{H_7}^{G}
\theta_7^{-k}$ are equal.  Thus $\Gamma_{H_7}$ is
\begin{eqnarray*}
\Gamma_{H_7} & = & \frac{1}{7} \left (\Ind_{H_7}^{G} \theta_7 + 2
\Ind_{H_7}^{G}
\theta_7^2 + \ldots + 6 \Ind_{H_7}^{G} \theta_7^6 \right ) \\
& = & \frac{1}{7} \left ( 7 \Ind_{H_7}^{G} \theta_7 + 7
\Ind_{H_7}^{G} \theta_7^2 + 7 \Ind_{H_7}^{G} \theta_7^4 \right ) \\
& = & \Ind_{H_7}^{G} \theta_7 + \Ind_{H_7}^{G} \theta_7^2 +
\Ind_{H_7}^{G} \theta_7^4.
\end{eqnarray*}

Recall from section \ref{sec:indchars} that the multiplicities of
the irreducible representations $W_\alpha$ and $X_\beta$ in the
induced representation $\Ind_{H_7}^G \theta_7^k$  depend on the
value of $\alpha(\phi)$ or $\beta(\phi)$, and that this value must
be $e^{\frac{2 \pi i k}{7}}$ for some $k = 0, \ldots, 6$. In the sum
$ \Gamma_{H_7} = \Ind_{H_7}^{G} \theta_7 + \Ind_{H_7}^{G} \theta_7^2
+ \Ind_{H_7}^{G} \theta_7^4$ we will have, for example for the
multiplicities of the $W_\alpha$ when $q \equiv 1 \pmod {28}$,

{\footnotesize{
\begin{eqnarray*}
\Gamma_{H_7} & = & \Ind_{H_7}^{G} \theta_7 + \Ind_{H_7}^{G}
\theta_7^2 + \Ind_{H_7}^{G} \theta_7^4 \\
& = & \frac{q-1}{7} \sum_{\alpha(\phi) \neq e^{\pm \frac{2 \pi i
}{7}}} W_{\alpha}+ \frac{q+6}{7} \sum_{\alpha(\phi) = e^{\pm \frac{2
\pi i }{7}}} W_{\alpha} \\
& + & \frac{q-1}{7} \sum_{\alpha(\phi) \neq e^{\pm \frac{4 \pi i
}{7}}} W_{\alpha} + \frac{q+6}{7} \sum_{\alpha(\phi) = e^{\pm
\frac{4 \pi i }{7}}} W_{\alpha} \\
&+& \frac{q-1}{7} \sum_{\alpha(\phi) \neq e^{\pm \frac{8 \pi i
}{7}}} W_{\alpha} + \frac{q+6}{7} \sum_{\alpha(\phi) = e^{\pm
\frac{8 \pi i }{7}}} W_{\alpha} \\
& + & \mbox{ other characters}.
\end{eqnarray*} }}
This adds up to
\begin{equation*}
\Gamma_{H_7} = \frac{3q+4}{7} \sum_{\alpha(\phi) \neq 1} W_{\alpha}+
\frac{3q-3}{7} \sum_{\alpha(\phi) = 1} W_{\alpha} + \mbox{ other
characters}.
\end{equation*}

The multiplicities of the other irreducible characters in
$\Ind_{H_7}^G \theta_7^k$ do not depend on $k$.  Adding these in,
the total for the case $q \equiv 1 \pmod {28}$ is

{\footnotesize{
\begin{eqnarray*}
\Gamma_{H_7} = \frac{3q-3}{7} \left [ \sum_{\beta} X_{\beta} + V +
\sum_{\alpha(\phi) =1} W_{\alpha} + \frac{1}{2} (W' + W'') \right ]
+ \frac{3q+4}{7} \sum_{\alpha(\phi) \neq 1} W_{\alpha}.
\end{eqnarray*} }}

Similar calculations yield the following.  If $q \equiv 13 \pmod
{28}$,

{\footnotesize{
\begin{eqnarray*}
\Gamma_{H_7} = \frac{3q+3}{7} \left [ \sum_{\beta(\phi) =1}
X_{\beta} + V + \sum_{\alpha} W_{\alpha} + \frac{1}{2} (W' + W'')
\right ] + \frac{3q-4}{7} \sum_{\beta(\phi) \neq 1} X_{\beta}.
\end{eqnarray*} }}

If $q \equiv 15 \pmod {28}$, we have

{\footnotesize{
\begin{eqnarray*}
\Gamma_{H_7} = \frac{3q-3}{7} \left [ \sum_{\beta} X_{\beta} + V +
\sum_{\alpha(\phi) =1} W_{\alpha} + \frac{1}{2} (X' + X'') \right ]
+ \frac{3q+4}{7} \sum_{\alpha(\phi) \neq 1} W_{\alpha}.
\end{eqnarray*} }}

And if $q \equiv 27 \pmod {28}$, we have

{\footnotesize{
\begin{eqnarray*}
\Gamma_{H_7} = \frac{3q+3}{7} \left [ \sum_{\beta(\phi) =1}
X_{\beta} + V + \sum_{\alpha} W_{\alpha} + \frac{1}{2} (X' + X'')
\right ] + \frac{3q-4}{7} \sum_{\beta(\phi) \neq 1} X_{\beta}.
\end{eqnarray*} }}

To compute the ramification module, we sum the components
$\Gamma_{H_2}$, $\Gamma_{H_3}$, and $\Gamma_{H_7}$ listed above. The
following numbers will be useful.

\begin{definition}
For each possible equivalence class of $q \pmod {84}$, we define a
\textbf{base multiplicity} $m$, as follows:
\begin{itemize}
\item  If $q \equiv 1, 13, 29, \mbox{ or } 43 \pmod {84}$, then $m =
q + \lfloor \frac{q}{84} \rfloor$.
\item  If $q \equiv 41, 55, 71, \mbox{ or } 83 \pmod {84}$, then $m =
q + \lceil \frac{q}{84} \rceil$.
\end{itemize}
\end{definition}

\begin{definition}
Let $\alpha: \fff^{\times} \to \ccc^{\times}$ be a character of
$\fff^{\times}$.  Then we define a number
\begin{equation*}
N_{\alpha} = \# \{ x \in \{i, \omega, \phi\} \ | \ x \in
\fff^{\times} \mbox{ and } \alpha(x) \neq 1\}.
\end{equation*}
\end{definition}

\begin{definition}
Recall that $T$ is the cyclic subgroup of
$\fff(\sqrt{\varepsilon})^{\times}$ of order $q+1$.  Let $\beta: T
\to \ccc^{\times}$ be a character of $T$. Then we define a number
\begin{equation*}
N_{\beta} = \# \{ x \in \{i, \omega, \phi\}\ |\ x \in T \mbox { and
} \beta(x) \neq 1\}.
\end{equation*}
\end{definition}

\begin{theorem}
\label{thrm:main} We have the following decomposition of the
ramification module:

\begin{itemize}
\item
If $q\equiv 1\pmod {8}$, then
\begin{equation*}
\tilde{\Gamma}_G = \frac{m}{2}(W'+W'')+ m V + \sum_{\beta} (m -
N_{\beta}) X_{\beta} + \sum_{\alpha} (m + N_{\alpha}) W_{\alpha}
\end{equation*}

\item
If $q\equiv 3\pmod {8}$, then
\begin{equation*}
\tilde{\Gamma}_G = \frac{m-1}{2}(X'+X'') + m V + \sum_{\beta}
(m-N_{\beta}) X_{\beta} + \sum_{\alpha} (m + N_{\alpha}) W_{\alpha}
\end{equation*}

\item
If $q\equiv 5\pmod {8}$, then
\begin{equation*}
\tilde{\Gamma}_G = \frac{m+1}{2}(W'+W'')+ m V + \sum_{\beta} (m -
N_{\beta}) X_{\beta} + \sum_{\alpha} (m + N_{\alpha}) W_{\alpha}
\end{equation*}

\item
If $q\equiv 7 \pmod {8}$, then
\begin{equation*}
\tilde{\Gamma}_G = \frac{m}{2}(X'+X'') + m V + \sum_{\beta}
(m-N_{\beta}) X_{\beta} + \sum_{\alpha} (m + N_{\alpha}) W_{\alpha}
\end{equation*}

\end{itemize}
\end{theorem}

\subsection{Equivariant degree}
\label{sec:equivdeg}

Now we will define and compute the equivariant degree of a
$G$-invariant divisor.  (See for example \cite{B} for more details).
This, together with the equivariant Riemann-Roch formula
(\ref{eqn:Borne}), will allow us to compute the $G$-module structure
of the Riemann-Roch space $L(D)$.

Fix a point $P\in X$ and let $D$ be a divisor on $X$ of the
form
\[
D=\frac{1}{e_P}\sum_{g\in G}g(P)=\sum_{g\in G/G_P}g(P),
\]
where $G_P$ denotes the stabilizer in $G$ of $P$ and $e_P=|G_P|$
denotes the ramification index at $P$. Such a divisor is called a
{\it reduced orbit}; any $G$-invariant divisor on $X$ can be written
as a sum of multiples of reduced orbits.

The {\it equivariant degree} of a multiple $rD$ of a reduced orbit
is the virtual representation
\begin{equation*}
\deg_{eq}(rD) = \left\{ \begin{array}{c c}  \Ind_{G_P}^G\
\sum_{\ell
=1}^r \theta_P^{-\ell}, & r > 0 \\
0, & r=0 \\
-\Ind_{G_P}^G\ \sum_{\ell =0}^{|r|-1} \theta_P^{-\ell}, & r < 0
\end{array} \right.
\end{equation*}
where $\theta_P$ is the ramification character of $X$ at $P$ (a
nontrivial character of $G_P$). In general, the equivariant degree
is additive on disjointly supported divisors.  Note that if $r$ is a
multiple of $e_P$, then then $D$ is the pull-back of a divisor on
$X/G$ via $\psi$ in (\ref{eqn:psi}), and the equivariant degree is a
multiple of the regular representation $\ccc[G]$ of $G$.  More
generally, if $D$ is a reduced orbit and $r = e_P r' + r ''$, then
\[
\deg_{eq}(rD) = r' \cdot \ccc[G] + \deg_{eq}(r''D).
\]
(Note this is true even when $r'$ is negative).

On the Hurwitz curve $X$, the results of section \ref{sec:ram} tell
us that there are only four types of reduced orbits to consider: the
stabilizer $G_P$ of a point $P$ in the support of $D$ may have order
$1$, $2$, $3$, or $7$, and therefore be either trivial or conjugate
to $H_2$, $H_3$, or $H_7$.  Let $D_1$, $D_2$, $D_3$, and $D_7$
denote reduced orbits of each type. There is only one choice of
reduced orbit for $D_2$, $D_3$, and $D_7$; for $D_1$ we see from the
definition that the equivariant degree does not depend on our choice
of orbit. Given a point in $D_1$, the stabilizer is trivial, so the
divisor is a pullback and the equivariant degree is
\[
\deg_{eq}(D_1) = \ccc[G].
\]

A general $G$-invariant divisor may be written as $r_1 D_1 + r_2 D_2
+ r_3 D_3 + r_7 D_7$.  If we write $r_2 = 2 r_2' + r_2 ''$, $r_3 = 3
r_3' + r_3 ''$, and $r_7 = 7 r_7' + r_7 ''$, then we have
\begin{equation*}
\begin{array}{l}
\deg_{eq}(r_1 D_1 + r_2 D_2 + r_3 D_3 + r_7 D_7) \\
\quad \quad = \deg_{eq}((r_1 + r_2' + r_3' + r_7') D_1 + r_2'' D_2 + r_3'' D_3 + r_7'' D_7)\\
\quad \quad = (r_1 + r_2' + r_3' + r_7')\ccc[G] + \deg_{eq}(r_2''
D_2 + r_3'' D_3 + r_7'' D_7).
\end{array}
\end{equation*}
Therefore, to compute the equivariant degree of a general divisor,
all that remains is to compute $\deg_{eq}(r_i D_i)$ for $i \in \{2,
3, 7\}$, where we may assume that $1 \leq r_i < i$.

\begin{description}

\item[Case 1]: $r_2 D_2$.
Given our assumptions, the only possibility is that $r_2=1$. Given a
point $P$ in the support of $D_2$, the stabilizer $G_P$ is conjugate
to $H_2$.  In this case, the equivariant degree of $D_2$ is
\[
\deg_{eq}(D_2) = \Ind_{H_2}^G \theta_2.
\]

\item[Case 2]: $r_3 D_3$.
Here we may have either $r_3 = 1$ or $r_3 = 2$.  The stabilizer of a
point in the support of $D_3$ is conjugate to $H_3$.  Recall that
$\Ind_{H_3}^G \theta_3^2= \Ind_{H_3}^G \theta_3$, so we have
\begin{eqnarray*}
\deg_{eq}(D_2) &=& \Ind_{H_3}^G \theta_3 \\
\deg_{eq}(2D_2)&=&  2 \Ind_{H_3}^G \theta_3.
\end{eqnarray*}

\item[Case 3]:  $r_7 D_7$.
In this case, we have $1 \leq r_7 \leq 6$.  The stabilizer of a
point in the support of $D_7$ is conjugate to $H_7$.  Recall that
for $k = 1, \ldots, 6$, $\Ind_{H_7}^G \theta_7^k = \Ind_{H_7}^G
\theta_7^{-k}$. Therefore the equivariant degree is as follows:
\begin{itemize}
\item $\deg_{eq}(D_7) = \Ind_{H_7}^G \theta_7$.
\item $\deg_{eq}(2 D_7) = \Ind_{H_7}^G \theta_7 + \Ind_{H_7}^G
\theta_7^2$.
\item $\deg_{eq}(3 D_7) = \Ind_{H_7}^G \theta_7 + \Ind_{H_7}^G
\theta_7^2 + \Ind_{H_7}^G \theta_7^3$, which is the same as the
$H_7$ component of the ramification module, $\Gamma_{H_7}$.
\item $\deg_{eq}(4 D_7) = \Gamma_{H_7} + \Ind_{H_7}^G \theta_7^3$.
\item $\deg_{eq}(5 D_7) = \Gamma_{H_7} + \Ind_{H_7}^G \theta_7^3 + \Ind_{H_7}^G \theta_7^2$.
\item $\deg_{eq}(6 D_7) = 2 \Gamma_{H_7}$.
\end{itemize}

\end{description}

Now we add these up.  As in the case of the ramification module, the
equivariant degree is most conveniently written in terms of a ``base
multiplicity'' and modifiers.  We define the base multiplicity as
follows.

\begin{itemize}
\item  If $q \equiv 1 \pmod 4$, then let $b_2 = r_2 \left (
\frac{q-1}{2} \right )$.  Otherwise, if $q \equiv 3 \pmod 4$, then
let $b_2 = r_2 \left ( \frac{q+1}{2} \right )$.

\item  If $q \equiv 1 \pmod 3$, then let $b_3 = r_3 \left (
\frac{q-1}{3} \right )$, and if $q \equiv 2 \pmod 3$, then let $b_3
= r_3 \left ( \frac{q+1}{3} \right )$.

\item  Similarly, if $q \equiv 1 \pmod 7$, then let $b_7 = r_7 \left (
\frac{q-1}{7} \right )$, and if $q \equiv 6 \pmod 7$, then let $b_7
= r_7 \left ( \frac{q+1}{7} \right )$.
\end{itemize}

The base multiplicity is then defined to be
\begin{eqnarray*}
b &=& b_2 + b_3 + b_7 \\
&=& r_2 \left ( \frac{q \pm 1}{2} \right ) + r_3 \left ( \frac{q \pm
1}{3} \right ) + r_7 \left ( \frac{q \pm 1}{7} \right ).
\end{eqnarray*}

Then the equivariant degree $\deg_{eq}(D)$ of the divisor $D=r_1 D_1
+ r_2 D_2 + r_3 D_3 + r_7 D_7$, with $0 \leq r_2 \leq 1$, $0 \leq
r_3 \leq 2$, and $0 \leq r_7 \leq 6$, is

\begin{equation}
\label{eqn:degeq}
\deg_{eq}(D) = b \left[ \sum_{\beta} X_{\beta} + V
+ \sum_{\alpha} W_{\alpha} \right] + \mbox{ modifiers},
\end{equation}
where the modifiers are listed in the table below.  For each $q$,
three of the rows below will be added.
\begin{center}
\begin{tabular}[ht]{|c|c|}
\hline $q$ & Modifiers to equivariant degree \\ \hline \hline
 $q \equiv 1 \pmod 8$   &  $ \displaystyle + \ r_2 \sum_{\alpha(i)=-1}
 W_{\alpha} + \frac{b}{2} (W'+W'')$ \\ \hline
$q \equiv 3 \pmod 8$   &  $ \displaystyle - \ r_2 \sum_{\beta(i)=-1}
 X_{\beta} + \frac{b-r_2}{2} (X'+X'')$ \\ \hline
 $q \equiv 5 \pmod 8$   &  $ \displaystyle + \ r_2 \sum_{\alpha(i)=-1}
 W_{\alpha} + \frac{b+r_2}{2} (W'+W'')$ \\ \hline
$q \equiv 7 \pmod 8$   &  $ \displaystyle - \ r_2 \sum_{\beta(i)=-1}
 X_{\beta} + \frac{b}{2} (X'+X'')$ \\ \hline \hline
$q \equiv 1 \pmod 3$   &  $ \displaystyle + \ r_3
\sum_{\alpha(\omega) \neq 1}  W_{\alpha}$ \\ \hline

$q \equiv 2 \pmod 3$   &  $ \displaystyle - r_3 \sum_{\beta(\omega)
\neq 1} X_{\beta} $ \\ \hline \hline

$q \equiv 1 \pmod 7$   &  $ \displaystyle + \ \sum_{k=1}^{r_7} \sum_{\alpha(\phi)= e^{\pm \frac{2 \pi i k}{7}} }  W_{\alpha}$ \\
\hline

$q \equiv 6 \pmod 7$   &  $ \displaystyle - \sum_{k=1}^{r_7}
\sum_{\beta(\phi) = e^{\pm \frac{2 \pi i k}{7}}} X_{\beta} $ \\
\hline \hline
\end{tabular}
\end{center}

\subsection{The Riemann-Roch space}

Now we would like to compute the $G$-module structure of the
Riemann-Roch space $L(D)$ for a $G$-invariant divisor $D$. First,
let us consider which $G$-invariant divisors are non-special. To be
non-special, it is sufficient to have $\deg D > 2g-2$, where

\[
g=1+\frac{(q)(q^2-1)}{168}
\]
is the genus of $X$, so $2g-2 =
\frac{1}{84}q(q^2-1)=\frac{1}{168}|G|$. The reduced orbits $D_1$,
$D_2$, $D_3$ and $D_7$ have degrees $|G|$, $|G|/2$, $|G|/3$, and
$|G|/7$, respectively.  Therefore if a $G$-invariant divisor $r_1
D_1 + r_2 D_2 + r_3 D_3 + r_7 D_7$ has positive degree, the smallest
its degree could be is $|G|/42$, which is strictly larger than
$2g-2$. Therefore any $G$-invariant divisor with positive degree is
non-special.

Thus for any $G$-invariant divisor $D$ with positive degree, we may
use the equivariant Riemann-Roch formula (\ref{eqn:Borne}) to
compute the $G$-module structure of the Riemann-Roch space $L(D)$:
\begin{equation*}
[L(D)]=(1-g_{X/G})[\ccc[G]]+[\deg_{eq}(D)]-[\tilde{\Gamma}_G].
\end{equation*}
Since $X/G \cong \ppp^1$, its genus is zero.  As in section
\ref{sec:equivdeg}, we may assume that $D=r_1 D_1 + r_2 D_2 + r_3
D_3 + r_7 D_7$, with $0 \leq r_2 \leq 1$, $0 \leq r_3 \leq 2$, and
$0 \leq r_7 \leq 6$.  Combining the results and notation of sections
\ref{sec:rammod} and \ref{sec:equivdeg}, we obtain the following.

\begin{equation*}
L(D) = (1+r_1) \ccc[G] + (b-m) \left[ \sum_{\beta} X_{\beta} + V +
\sum_{\alpha} W_{\alpha} \right] + \mbox{ modifiers},
\end{equation*}
where the modifiers depend on $q \pmod {168}$ and are listed in the
following table.  Again, for each value of $q$, three of the rows
below will be added.

\begin{center}
\begin{tabular}[ht]{|c|c|}
\hline $q$ & Modifiers to Riemann-Roch space \\ \hline \hline
 $q \equiv 1 \pmod 8$   &  $ \displaystyle + \ (r_2-1) \sum_{\alpha(i)=-1}
 W_{\alpha} + \frac{b-m}{2} (W'+W'')$ \\ \hline
$q \equiv 3 \pmod 8$   &  $ \displaystyle + \ (1-r_2)
\sum_{\beta(i)=-1}  X_{\beta} + \frac{b-m + 1 -r_2}{2} (X'+X'')$ \\
\hline
 $q \equiv 5 \pmod 8$   &  $ \displaystyle + \ (r_2-1) \sum_{\alpha(i)=-1}
 W_{\alpha} + \frac{b-m+r_2-1}{2} (W'+W'')$ \\ \hline
$q \equiv 7 \pmod 8$   &  $ \displaystyle + \ (1-r_2)
\sum_{\beta(i)=-1}
 X_{\beta} + \frac{b-m}{2} (X'+X'')$ \\ \hline \hline
$q \equiv 1 \pmod 3$   &  $ \displaystyle + \ (r_3-1)
\sum_{\alpha(\omega) \neq 1}  W_{\alpha}$ \\ \hline

$q \equiv 2 \pmod 3$   &  $ \displaystyle + \ (1- r_3)
\sum_{\beta(\omega) \neq 1} X_{\beta} $ \\ \hline \hline

$q \equiv 1 \pmod 7$   &  $ \displaystyle + \ \sum_{k=1}^{r_7}
\sum_{\alpha(\phi)= e^{\pm \frac{2 \pi i k}{7}} }  W_{\alpha} -
\sum_{alpha(\phi) \neq 1} W_{\alpha} $ \\ \hline

$q \equiv 6 \pmod 7$   &  $ \displaystyle + \ \sum_{\beta(\phi) \neq
1} X_{\beta} - \sum_{k=1}^{r_7}
\sum_{\beta(\phi) = e^{\pm \frac{2 \pi i k}{7}}} X_{\beta} $ \\
\hline \hline
\end{tabular}
\end{center}

\subsection{Action on holomorphic differentials}
\label{sec:cohomology}
As a corollary, it is an easy exercise now to
compute explicitly the decomposition

$$
H^1(X,\ccc) = H^0(X,\Omega^1)\oplus \overline{H^0(X,\Omega^1)}
=L(K_X)\oplus \overline{L(K_X)},
$$
into irreducible $G$-modules, where $K_X$ is a canonical divisor of
$X$. The action of $G$ on the complex conjugate vector space
$\overline{L(K_X)}$ of $L(K_X)$ will be by the complex conjugate
(contragredient) representation.  The Riemann-Hurwitz theorem tells
us that
\begin{eqnarray*}
K_X &=& \pi^*(K_{\ppp^1}) + R\\
    &=& -2 D_1 + D_2 + 2 D_3 + 6 D_7
\end{eqnarray*}
where $R$ is the ramification divisor.  Thus the equivariant degree
of $K_X$ is $\deg_{eq}(K_X) = -2 \cdot \ccc[G] + \deg_{eq}(R)$.
Note from the preliminary equivariant degree calculations, that
\begin{eqnarray*}
\deg_{eq}(R) &=& \deg_{eq}{D_2} + \deg_{eq}{2 D_3} + \deg_{eq}{6
D_7} \\
&=& \Ind_{H_2}^{G} \theta_2 + 2 \Ind_{H_3}^{G} \theta_3 +
\sum_{k=1}^6 \Ind_{H_7}^{G} \theta_7^k \\
&=& 2 \Gamma_{H_2} + 2 \Gamma_{H_3} + 2 \Gamma_{H_7} \\
&=& 2 \tilde{\Gamma}.
\end{eqnarray*}

Therefore, using the equivariant Riemann-Roch formula
(\ref{eqn:Borne}),
\begin{equation}
L(K_X) = \tilde{\Gamma} - \ccc[G].
\end{equation}
We will see in the next section that this is invariant under complex
conjugation, so that as $G$-modules, $H^1(X,\ccc) \cong 2L(K_X)$.

Using the results of section \ref{sec:rammod}, we obtain the
following.

\begin{theorem}
The $G$-module structure of $L(K)=H^0(X,\Omega^1)$ is as follows:
\begin{itemize}
\item If $q \equiv 1$, $97$, or $113 \pmod {168}$, then
{\footnotesize{
\begin{equation*}
L(K_X) = \frac{\lfloor \frac{q}{84} \rfloor -1}{2} (W' + W'') +
\sum_{\beta} \left( \lfloor \frac{q}{84} \rfloor + 1 - N_{\beta}
\right) X_{\beta} + \lfloor \frac{q}{84} \rfloor V  + \sum_{\alpha}
\left( \lfloor \frac{q}{84} \rfloor - 1 + N_{\alpha} \right)
W_{\alpha}.
\end{equation*} }}

\item If $q \equiv 43 \pmod {168}$, then
{\footnotesize{ \begin{equation*} L(K_X) = \lfloor \frac{q}{84}
\rfloor \left[ \frac{1}{2}(X' + X'') + V \right ] + \sum_{\beta}
\left( \lfloor \frac{q}{84} \rfloor + 1 - N_{\beta} \right)
X_{\beta} + \sum_{\alpha} \left( \lfloor \frac{q}{84} \rfloor - 1 +
N_{\alpha} \right) W_{\alpha}.
\end{equation*} }}

\item If $q \equiv 13$, $29$, or $85 \pmod {168}$, then
{\footnotesize{ \begin{equation*} L(K_X) = \lfloor \frac{q}{84}
\rfloor \left[ \frac{1}{2}(W' + W'') + V \right ] + \sum_{\beta}
\left( \lfloor \frac{q}{84} \rfloor + 1 - N_{\beta} \right)
X_{\beta} + \sum_{\alpha} \left( \lfloor \frac{q}{84} \rfloor - 1 +
N_{\alpha} \right) W_{\alpha}.
\end{equation*} }}

\item If $q \equiv 127 \pmod {168}$, then
{\footnotesize{ \begin{equation*} L(K_X) = \frac{\lfloor
\frac{q}{84} \rfloor +1}{2} (X' + X'') + \sum_{\beta} \left( \lfloor
\frac{q}{84} \rfloor + 1 - N_{\beta} \right) X_{\beta} + \lfloor
\frac{q}{84} \rfloor V + \sum_{\alpha} \left( \lfloor \frac{q}{84}
\rfloor - 1 + N_{\alpha} \right) W_{\alpha}.
\end{equation*} }}

\item If $q \equiv 41 \pmod {168}$, then
{\footnotesize{ \begin{equation*} L(K_X) = \frac{\lceil \frac{q}{84}
\rceil -1}{2} (W' + W'') + \sum_{\beta} \left( \lceil \frac{q}{84}
\rceil + 1 - N_{\beta} \right) X_{\beta} + \lceil \frac{q}{84}
\rceil V  + \sum_{\alpha} \left( \lceil \frac{q}{84} \rceil - 1 +
N_{\alpha} \right) W_{\alpha}.
\end{equation*} }}

\item If $q \equiv 83$, $139$, or $155 \pmod {168}$, then
{\footnotesize{ \begin{equation*} L(K_X) = \lceil \frac{q}{84}
\rceil \left[ \frac{1}{2}(X' + X'') + V \right ] + \sum_{\beta}
\left( \lceil \frac{q}{84} \rceil + 1 - N_{\beta} \right) X_{\beta}
+ \sum_{\alpha} \left( \lceil \frac{q}{84} \rceil - 1 + N_{\alpha}
\right) W_{\alpha}.
\end{equation*} }}

\item If $q \equiv 125 \pmod {168}$, then
{\footnotesize{ \begin{equation*} L(K_X) = \lceil \frac{q}{84}
\rceil \left[ \frac{1}{2}(W' + W'') + V \right ] + \sum_{\beta}
\left( \lceil \frac{q}{84} \rceil + 1 - N_{\beta} \right) X_{\beta}
+ \sum_{\alpha} \left( \lceil \frac{q}{84} \rceil - 1 + N_{\alpha}
\right) W_{\alpha}.
\end{equation*} }}

\item If $q \equiv 55$, $71$, or $167$ $\pmod {168}$, then
{\footnotesize{ \begin{equation*} L(K_X) = \frac{\lceil \frac{q}{84}
\rceil +1}{2} (X' + X'') + \sum_{\beta} \left( \lceil \frac{q}{84}
\rceil + 1 - N_{\beta} \right) X_{\beta} + \lceil \frac{q}{84}
\rceil V  + \sum_{\alpha} \left( \lceil \frac{q}{84} \rceil - 1 +
N_{\alpha} \right) W_{\alpha}.
\end{equation*} }}

\end{itemize}
\end{theorem}

\section{Galois action}

As discussed in section \ref{sec:rep thy}, there is a Galois action
on the set of equivalence classes of irreducible representations of
$PSL(2,q)$.  One question of obvious interest is whether the modules
we have computed are invariant under this action.

\begin{theorem}  The ramification module is Galois-invariant.
\end{theorem}

\pf  Recall from section \ref{sec:rep thy} that the Galois group
$\mathcal{G}$ permutes $m$th roots of unity, where $m = q(q^2-1)/4$.
It acts on representations of $PSL(2,q)$ by permuting character
values.  Thus it fixes the trivial representation and the
$q$-dimensional representation $V$, whose character values are
rational.  It will act as a permutation on the representations
$W_{\alpha}$ and on the representations $X_{\beta}$.  Lastly, it
will act as an involution on either the representations $W'$ and
$W''$ or $X'$ and $X''$.

Because the multiplicities of $W'$ and $W''$ or $X'$ and $X''$ are
the same in the ramification module, the Galois action will be
invariant on this component.  The multiplicity of a representation
$W_{\alpha}$ or $X_{\beta}$ in the ramification module depends on
the number $N_{\alpha}$ or $N_{\beta}$, which is determined by the
value of the character $\alpha$ or $\beta$ on the special numbers
$i$, $\omega$, and $\phi$.  In fact, the numbers $N_{\alpha}$ and
$N_{\beta}$ are determined only by whether these character values
are equal to $1$ or not equal to $1$.  Since an element of the
Galois group will take a character value to a power of itself, the
Galois action must preserve the numbers $N_{\alpha}$ and
$N_{\beta}$.  Therefore this component of the ramification module is
invariant as well. \qed

Since the ramification module is Galois-invariant, and of course the
regular representation is Galois-invariant, $L(K_X)$ will be Galois
invariant.  In particular, as stated in section
\ref{sec:cohomology}, $L(K_X)$ will be invariant under complex
conjugation. For a general divisor $D$, the Riemann-Roch space
$L(D)$ will be Galois-invariant if and only if the equivariant
degree of $D$ is.

\begin{theorem}
Let $D=r_1 D_1 + r_2 D_2 + r_3 D_3 + r_7 D_7$ be a $G$-invariant
divisor.  Then the equivariant degree of $D$ is Galois-invariant if
$r_7 \in \{0, 3, 6\} \pmod{7}$.
\end{theorem}

\pf   As in section \ref{sec:equivdeg}, multiples of $2$ in $r_2$,
$3$ in $r_3$, and $7$ in $r_7$ can be absorbed into the $r_1 D_1$
term without affecting the equivariant degree.  Therefore we may
assume that $0 \leq r_2 \leq 1$, $0 \leq r_3 \leq 2$, and $0 \leq
r_7 \leq 6$.

The result can again be seen by looking at the multiplicities of
representations permuted by the Galois group. The multiplicities of
$W'$ and $W''$ or $X'$ and $X''$ are the same. By (\ref{eqn:degeq}),
the multiplicity of a representation $W_{\alpha}$ or $X_{\beta}$
depends on $r_2$, $r_3$, and $r_7$, and not on $r_1$.  Again, the
Galois action will not permute a representation $W_{\alpha}$ with
$\alpha(i)=1$ with one with $\alpha(i) \neq 1$; similarly for
$X_{\beta}$, and for $\omega$. However, it could permute for example
a representation $W_{\alpha}$ with $\alpha(\phi)=e^{\frac{2 \pi
i}{7}}$ with one with $\alpha(\phi)=e^{\frac{4 \pi i}{7}}$.  Thus
the equivariant degree may not be Galois-invariant unless the
multiplicities of these representations are equal.  In the cases
where $r_7 \in \{0, 3, 6\}$, then these multiplicities will be
equal; otherwise they will not. \qed

Note that for some values of $q$, the equivariant degree may be
Galois-invariant even if $r_7$ is not $0$, $3$, or $6$.

A previous result of the first two authors (see \cite{JK1}) gives a
simpler formula (see equation \ref{eqn:JKrammod}) to compute the
multiplicity of an irreducible representation in the ramification
module, when the ramification module is Galois-invariant. In the
example at hand, if $r_7 \in \{0, 3, 6\}$, then since the
equivariant degree is a multiple of the $H_7$ component of the
ramification module, a slight modification of this formula gives an
easy computation of the equivariant degree and therefore the
Riemann-Roch space.

\begin{corollary}
Let $D = r_1 D_1 + r_2 D_2 + r_3 D_3 + r_7 D_7$, with $0 \leq r_2
\leq 1$, $0 \leq r_3 \leq 2$, and $r_7 \in \{0, 3, 6\}$.  Then
\begin{eqnarray*}
L(D) &=& \bigoplus_{\pi \in G^*} \left[ (1+ r_1 + r_2 +
\frac{r_3}{2} +
\frac{r_7}{6}) \dim \pi  \right.\\
&& \left. + (\frac{1}{2}- r_2) \dim \pi^{H_2} + (\frac{1}{2}-
\frac{r_3}{2}) \dim \pi^{H_3} +(\frac{1}{2}- \frac{r_7}{6}) \dim
\pi^{H_7} \right] \pi.
\end{eqnarray*}
\end{corollary}

Note that in spite of appearances, the multiplicity of each
irreducible representation will in fact be an integer.

\pf We see from the calculations in section \ref{sec:equivdeg} that
the equivariant degree of $D$ is equal to
\begin{eqnarray*}
\deg_{eq}(D) &=& r_1 \ccc[G] + 2 r_2 \Gamma_{H_2} + r_3 \Gamma_{H_3}
+
\frac{r_7}{3} \Gamma_{H_7} \\
&=& \bigoplus_{\pi \in G^*} \left[ (r_1 + r_2 + \frac{r_3}{2} +
\frac{r_7}{6}) \dim \pi  \right.\\
&& \left. - r_2 \dim \pi^{H_2} - \frac{r_3}{2} \dim \pi^{H_3} -
\frac{r_7}{6} \dim \pi^{H_7} \right] \pi.
\end{eqnarray*}

The ramification module is
\begin{equation*}
\tilde{\Gamma}_G =\bigoplus_{\pi\in G^*} \left[ \sum_{\ell \in
{2,3,7}} ({\rm dim}\, \pi -{\rm dim}\, (\pi^{H_\ell}))
\frac{1}{2}\right]\pi .
\end{equation*}
This sum splits into $\tilde{\Gamma}_G = \Gamma_{H_2} + \Gamma_{H_3}
+ \Gamma_{H_7}$ in the obvious way along the inner sum. Putting
these together using the equivariant Riemann-Roch formula
(\ref{eqn:Borne}), we obtain the desired result. \qed

\end{document}